\documentclass[a4paper,10 pt]{article}
\usepackage{amsmath}
\usepackage{amssymb}
\usepackage{amsthm}
\usepackage{enumerate}
\usepackage{amsfonts}
\usepackage{mathrsfs}
\begin{document}
\title{\Large\textbf{A NOTE ABOUT THE $(2,3)$-GENERATION OF $SL_{12}(q)$}}
\author{TSANKO RAYKOV GENCHEV}
\date{}
\maketitle
\begin{abstract}
In this note we provide  explicit uniform type $(2,3)$-generators for the  special linear group $SL_{12}(q)$ for all $q$ except for $q=2$ or $q=4$. Our considerations are easily traceable, self-contained and based only on the known list of maximal subgroups of this group.
\end{abstract}
\vspace{16pt}
\indent\indent Now, when the problem concerning the $(2,3)$-generation (especially) of the special linear groups and their projective images is completely solved (see \cite{4} and references there), we give our contribution by discussing the last remaining group $SL_{12}(q)$ in the light of the method used in our works \cite{2}, \cite{3} and \cite{5}. It has been proved in \cite{4} the $(2,3)$-generation of the same group over all finite fields in characteristic not equal to $5$. The author's approach is different from ours in which we make an essential use of the known list of maximal subgroups of $SL_{12}(q)$.\\
\indent The group $G = SL_{12}(q)$, where $q = p^{m}$ and $p$ is a prime, acts naturally on a twelve-dimensional vector space $V$ over the field $F = GF(q)$. We identify $V$ with the column vectors of $F^{12}$, and let $v_{1}$, . . . , $v_{12}$ be the standard base of the space $V$. Set $Q = q^{11}-1$ if $q \neq 3, 7$ and $Q = (q^{11}-1)/2$ if  $q = 3, 7$.\\
\indent First of all, based only to the orders of the maximal subgroups of $G$ (which exact structure is given in Tables $8.76$ and $8.77$ in \cite{1}), it can be easily seen that there is only one class of such subgroups having  an element of order $Q$. Namely, this is the class of reducible on the space $V$ subgroups; actually they are the stabilizers in $G$ of one or eleven-dimensional subspaces of $V$.\\
\indent Now, let us choose an element $\omega$ of order $Q$ in the multiplicative group of the field $GF(q^{11})$ and set 
\begin{center}
$f(t) = (t - \omega)(t - \omega^{q})(t - \omega^{q^{2}})(t - \omega^{q^{3}})(t - \omega^{q^{4}})(t - \omega^{q^{5}})(t - \omega^{q^{6}})(t - \omega^{q^{7}})(t - \omega^{q^{8}})(t - \omega^{q^{9}}) (t - \omega^{q^{10}}) = t^{11} - \alpha_{1}t^{10} + \alpha_{2}t^{9} - \alpha_{3}t^{8} + \alpha_{4}t^{7} - \alpha_{5}t^{6} + \alpha_{6}t^{5} -\alpha_{7}t^{4} + \alpha_{8}t^{3} - \alpha_{9}t^{2} + \alpha_{10}t - \alpha_{11}$.
\end{center}
Then $f(t) \in F[t]$ and the polynomial $f(t)$ is irreducible over the field $F$. Note that  $\alpha_{11} = \omega^\frac{q^{11} - 1}{q - 1}$  has order  $q - 1$ if $q \neq 3, 7$, $\alpha_{11} = 1$ if $q = 3$, and  $\alpha_{11}^{3} = 1 \neq \alpha_{11}$  if  $q = 7$.\\
\indent Set
\[ x = \left[ \begin{array}{cccccccccccc}
-1 & 0 & 0 & 0 & 0 & 0 & 0 & 0 & 0 & \alpha_{6}\alpha_{11}^{-1} & 0 & \alpha_{6}\\
0 & -1 & 0 & 0 & 0 & 0 & 0 & 0 & 0 & \alpha_{5}\alpha_{11}^{-1} & 0 & \alpha_{5}\\
0 & 0 & 0 & 0 & 0 & -1 & 0 & 0 & 0 & \alpha_{4}\alpha_{11}^{-1} & 0 & \alpha_{7}\\
0 & 0 & 0 & 0 & 0 & 0 & 0 & -1 & 0 & \alpha_{3}\alpha_{11}^{-1} & 0 & \alpha_{9}\\
0 & 0 & 0 & 0 & -1 & 0 & 0 & 0 & 0 & \alpha_{8}\alpha_{11}^{-1} & 0 & \alpha_{8}\\
0 & 0 & -1 & 0 & 0 & 0 & 0 & 0 & 0 & \alpha_{7}\alpha_{11}^{-1} & 0 & \alpha_{4}\\
0 & 0 & 0 & 0 & 0 & 0 & 0 & 0 & 0 & \alpha_{1}\alpha_{11}^{-1} & -1 & \alpha_{10}\\
0 & 0 & 0 & -1 & 0 & 0 & 0 & 0 & 0 & \alpha_{9}\alpha_{11}^{-1} & 0 & \alpha_{3}\\
0 & 0 & 0 & 0 & 0 & 0 & 0 & 0 & -1 & \alpha_{2}\alpha_{11}^{-1} & 0 & \alpha_{2}\\
0 & 0 & 0 & 0 & 0 & 0 & 0 & 0 & 0 & 0 & 0 & \alpha_{11}\\
0 & 0 & 0 & 0 & 0 & 0 & -1 & 0 & 0 & \alpha_{10}\alpha_{11}^{-1} & 0 & \alpha_{1}\\
0 & 0 & 0 & 0 & 0 & 0 & 0 & 0 & 0 & \alpha_{11}^{-1} & 0 & 0\\
 \end{array} \right],\] 
\[y = \left[ \begin{array}{cccccccccccc}
0 & 0 & 1 & 0 & 0 & 0 & 0 & 0 & 0 & 0 & 0 & 0\\
1 & 0 & 0 & 0 & 0 & 0 & 0 & 0 & 0 & 0 & 0 & 0\\
0 & 1 & 0 & 0 & 0 & 0 & 0 & 0 & 0 & 0 & 0 & 0\\
0 & 0 & 0 & 0 & 0 & 1 & 0 & 0 & 0 & 0 & 0 & 0\\
0 & 0 & 0 & 1 & 0 & 0 & 0 & 0 & 0 & 0 & 0 & 0\\
0 & 0 & 0 & 0 & 1 & 0 & 0 & 0 & 0 & 0 & 0 & 0\\
0 & 0 & 0 & 0 & 0 & 0 & 0 & 0 & 1 & 0 & 0 & 0\\
0 & 0 & 0 & 0 & 0 & 0 & 1 & 0 & 0 & 0 & 0 & 0\\
0 & 0 & 0 & 0 & 0 & 0 & 0 & 1 & 0 & 0 & 0 & 0\\
0 & 0 & 0 & 0 & 0 & 0 & 0 & 0 & 0 & 0 & 0 & 1\\
0 & 0 & 0 & 0 & 0 & 0 & 0 & 0 & 0 & 1 & 0 & 0\\
0 & 0 & 0 & 0 & 0 & 0 & 0 & 0 & 0 & 0 & 1 & 0\\
\end{array} \right].\]
Then $x$ and $y$ are elements of $G$ of orders $2$ and $3$, respectively. Denote 
\[z = xy = \left[ \begin{array}{cccccccccccc}
0 & 0 & -1 & 0 & 0 & 0 & 0 & 0 & 0 & 0 & \alpha_{6} & \alpha_{6}\alpha_{11}^{-1}\\
-1 & 0 & 0 & 0 & 0 & 0 & 0 & 0 & 0 & 0 & \alpha_{5} & \alpha_{5}\alpha_{11}^{-1}\\
0 & 0 & 0 & 0 & -1 & 0 & 0 & 0 & 0 & 0 & \alpha_{7} & \alpha_{4}\alpha_{11}^{-1}\\
0 & 0 & 0 & 0 & 0 & 0 & -1 & 0 & 0 & 0 & \alpha_{9} & \alpha_{3}\alpha_{11}^{-1}\\
0 & 0 & 0 & -1 & 0 & 0 & 0 & 0 & 0 & 0 & \alpha_{8} & \alpha_{8}\alpha_{11}^{-1}\\
0 & -1 & 0 & 0 & 0 & 0 & 0 & 0 & 0 & 0 & \alpha_{4} & \alpha_{7}\alpha_{11}^{-1}\\
0 & 0 & 0 & 0 & 0 & 0 & 0 & 0 & 0 & -1 & \alpha_{10} & \alpha_{1}\alpha_{11}^{-1}\\
0 & 0 & 0 & 0 & 0 & -1 & 0 & 0 & 0 & 0 & \alpha_{3} & \alpha_{9}\alpha_{11}^{-1}\\
0 & 0 & 0 & 0 & 0 & 0 & 0 & -1 & 0 & 0 & \alpha_{2} & \alpha_{2}\alpha_{11}^{-1}\\
0 & 0 & 0 & 0 & 0 & 0 & 0 & 0 & 0 & 0 & \alpha_{11} & 0\\
0 & 0 & 0 & 0 & 0 & 0 & 0 & 0 & -1 & 0 & \alpha_{1} & \alpha_{10}\alpha_{11}^{-1}\\
0 & 0 & 0 & 0 & 0 & 0 & 0 & 0 & 0 & 0 & 0 & \alpha_{11}^{-1}\\ 
\end{array} \right].\]
The characteristic polynomial of $z$ is $f_ {z}(t) = (t - \alpha_{11}^{-1})f(t)$ and the characteristic roots $\alpha_{11}^{-1}$, $\omega$, $\omega^{q}$, $\omega^{q^{2}}$, $\omega^{q^{3}}$, $\omega^{q^{4}}$, $\omega^{q^{5}}$, $\omega^{q^{6}}$, $\omega^{q^{7}}$, $\omega^{q^{8}}$, $\omega^{q^{9}}$, and $\omega^{q^{10}}$ of $z$ are pairwise distinct. Then, in $GL_{12}(q^{11})$, $z$ is conjugate to the matrix diag ($\alpha_{11}^{-1}$, $\omega$, $\omega^{q}$, $\omega^{q^{2}}$, $\omega^{q^{3}}$, $\omega^{q^{4}}$, $\omega^{q^{5}}$, $\omega^{q^{6}}$, $\omega^{q^{7}}$, $\omega^{q^{8}}$, $\omega^{q^{9}}$, $\omega^{q^{10}}$) and hence $z$ is an element of $SL_{12}(q)$ of order $Q$.\\
\indent Let $H = \left\langle {x,y}\right\rangle$, $H \leq G$. We prove that $H$ acts irreducibly on the space $V$. Indeed, assume that $W$ is an $H$-invariant subspace of $V$ and $k$ = dim $W$, $1 \leq k \leq 11$.\\
\indent Let first $k = 1$ and $0 \neq w \in W$. Then $y(w) = \lambda w$ where $\lambda \in F$ and $\lambda^{3} = 1$. This yields
\begin{center}
$w = \mu_{1}(v_{1} + \lambda^{2} v_{2} + \lambda v_{3}) + \mu_{2}(v_{4} + \lambda^{2} v_{5} + \lambda v_{6}) + \mu_{3}(v_{7} + \lambda^{2} v_{8} + \lambda v_{9}) + \mu_{4}(v_{10} + \lambda^{2} v_{11} + \lambda v_{12})$,  
\end{center}
where $\mu_{1}, \mu_{2}, \mu_{3},$ and $\mu_{4}$ are elements of the field $F$.\\
Now, we involve the action of $x$ onto $w$: $x(w) = \nu w$ where $\nu = \pm 1$. This yields consecutively $\mu_{4} \neq 0$ , $\alpha_{11} = \lambda^{2}\nu$, and 
\begin{enumerate}[\label=(1)]
  \item \rule{0pt}{0pt}\vspace*{-12pt}
    \begin{equation*}
     \mu_{3} = \lambda(\alpha_{1} + \nu \alpha_{10} -\lambda\nu)\mu_{4},
    \end{equation*}
\end{enumerate}
\begin{enumerate}[\label =(2)]
  \item \rule{0pt}{0pt}\vspace*{-12pt}
    \begin{equation*}
     \nu\mu_{1} + \mu_{2} = (\nu\alpha_{4}+ \alpha_{7})\mu_{4},
    \end{equation*}
\end{enumerate}
\begin{enumerate}[\label =(3)]
  \item \rule{0pt}{0pt}\vspace*{-12pt}
    \begin{equation*}
     \nu\mu_{2} + \lambda^{2}\mu_{3} = \lambda(\nu\alpha_{3}+ \alpha_{9})\mu_{4},
    \end{equation*}
\end{enumerate}
\begin{enumerate}[\label =(4)]
  \item \rule{0pt}{0pt}\vspace*{-12pt}
    \begin{equation*}
     (\nu + 1)(\mu_{1} - \lambda\alpha _{6}\mu_{4}) = 0,
    \end{equation*}
\end{enumerate}
\begin{enumerate}[\label =(5)]
  \item \rule{0pt}{0pt}\vspace*{-12pt}
    \begin{equation*}
      (\nu + 1)(\mu_{1} - \lambda^{2}\alpha _{5}\mu_{4}) = 0,
    \end{equation*}
\end{enumerate}
\begin{enumerate}[\label =(6)]
  \item \rule{0pt}{0pt}\vspace*{-12pt}
    \begin{equation*}
      (\nu + 1)(\mu_{2} - \lambda^{2}\alpha _{8}\mu_{4}) = 0,
    \end{equation*}
\end{enumerate}
\begin{enumerate}[\label =(7)]
  \item \rule{0pt}{0pt}\vspace*{-12pt}
    \begin{equation*}
      (\nu + 1)(\mu_{3} - \alpha _{2}\mu_{4}) = 0,
    \end{equation*}
\end{enumerate}
\noindent In particular, we have $\alpha_{11}^{3} = \nu$ and $\alpha_{11}^{6} = 1$. This is impossible if $q = 5$ or $q > 7$ since then $\alpha_{11}$ has order $q - 1$.
 According to our assumption ($q \neq 2, 4$) only two possibilities left: $q = 3$ (and $\alpha_{11} = 1$) or $q = 7$ (and $\alpha_{11}^{3} = 1 \neq \alpha_{11}$). So $\nu = 1$, $\alpha_{11} = \lambda^{2}$ and (1), (2), (3), (4), (5), (6) and (7) produce $\alpha_{1} = \lambda^{2}\alpha_{2} - \alpha_{10} + \lambda$, $\alpha_{6} = \lambda\alpha_{5}$, $\alpha_{7} = \lambda^{2}\alpha_{5} + \lambda^{2}\alpha_{8} - \alpha_{4}$ and $\alpha_{9} = \lambda\alpha_{2} + \lambda\alpha_{8} - \alpha_{3}$. Now $f(-1) =  -(1 + \lambda + \lambda^{2})(1 + \alpha_{2} + \alpha_{5} + \alpha_{8}) = 0$ both for $q = 3$ and $q = 7$, an impossibility as $f(t)$ is irreducible over the field $F$.\\
\indent Now let $2 \leq k \leq 11$. Then the characteristic polynomial of $z|_{W}$ has degree $k$ and has to divide $f_ {z}(t) = (t - \alpha_{11}^{-1})f(t)$.The irreducibility of $f(t)$ over $F$ leads immediately to the conclusion that this polynomial is $f(t)$ and $k=11$. Now the subspace $U$ of $V$ which is generated by the vectors $v_{1}$, $v_{2}$, $v_{3}$, . . .  , $v_{11}$ is $\left\langle {z}\right\rangle$-invariant. If $W \neq U$ then $U \cap W$ is $\left\langle {z}\right\rangle$-invariant and dim $(U \cap W) = 10$. This means that the characteristic polynomial of $z|_{U \cap W}$ has degree $10$ and must divide $f_ {z}(t) = (t - \alpha_{11}^{-1})f(t)$ which is impossible. Thus $W = U$ but obviously $U$ is not $\left\langle {y}\right\rangle$-invariant, a contradiction.\\
\indent Note that the above considerations fail if $q = 2$ or $4$.\\
\indent Now, as $H = \left\langle{x,y}\right\rangle$ acts irreducible on the space $V$ and it has an element of order $Q$, we conclude that $H$ can not be contained in any maximal subgroup of $G (= SL_{12}(q))$. Thus $H = G$ and $G = \left\langle {x,y}\right\rangle$ is a $(2,3)$-generated group. \hfill $\square$

\begin{center}

\end{center}
Ts. Genchev\\ 
e-mail: genchev57@yahoo.com\\
Department of Mathematics\\                                                              
Technical University\\
Varna, Bulgaria

\begin{thebibliography}{99}


\bibitem{1}
J. N. BRAY, D.F. HOLT, C. M. RONEY-DOUGAL. The Maximal Subgroups of the Low - Dimensional Finite Classical Groups. London Math. Soc. Lecture Note Series \textbf {407}, Cambridge University Press (2013).

\bibitem{2}
TS. GENCHEV, E. GENCHEVA. $(2,3)$-generation of the special linear groups of dimension $8$. \emph{Proceedings of the Forty Fourth Spring Conference of the Union of Bulgarian Mathematicians, SOK "Kamchia"}, April 2-6 (2015), 167-173.

\bibitem{3}
E. GENCHEVA, TS. GENCHEV and K. TABAKOV. $(2,3)$-generation of the special linear groups of dimensions $9$, $10$ and $11$. \emph{arXiv} \textbf{1412.8631v5} (2016).

\bibitem {4} 
M. A. PELLEGRINI. The $(2,3)$-generation of the special linear groups over finite fields. \emph{arXiv} \textbf{1605.04276v1} (2016).

\bibitem{5}
K. TABAKOV, E. GENCHEVA and TS. GENCHEV. $(2,3)$-generation of the special linear groups of dimension $11$. \emph{Proceedings of the Forty Fifth Spring Conference of the Union of Bulgarian Mathematicians, Pleven}, April 6-10 (2016), 146-151.



\end{thebibliography}
\end{document}